\documentclass[11pt]{amsart}

\usepackage{amsmath,amsfonts,amssymb,amscd}
\begin{document}
\renewcommand{\thefootnote}{\fnsymbol{footnote}}
\pagestyle{plain}

\title{Asymptotic polybalanced kernels\\
on extremal K\"ahler manifolds}
\author{Toshiki Mabuchi${}^*$}
\maketitle
\footnotetext{ ${}^{*}$Supported 
by JSPS Grant-in-Aid for Scientific Research (B) No. 25287010.}
\centerline{\it In honor of Professor Ngaiming Mok's 60th birthday}
\abstract
In this paper, improving a result 
 in \cite{M0}, we obtain asymptotic polybalanced 
kernels associated to extremal K\"ahler metrics on polarized algebraic manifolds.
As a corollary, 
we strengthen a result in \cite{M3} on 
asymptotic relative Chow-polystability for extremal K\"ahler polarized
algebraic manifolds. Finally, related to the Yau-Tian-Donaldson Conjecture for extremal K\"ahler metrics,
we shall discuss the difference between strong relative K-stability (cf.~\cite{MN}) and relative K-stability.
\endabstract

\section{Introduction}

In this paper, we fix once for all a {\it polarized algebraic manifold} $(M,L)$ 
which is by definition a pair of a nonsingular irreducible complex projective variety $M$ of dimension $n$
and a very ample holomorphic line bundle $L$ on $M$.
By taking the identity component $\operatorname{Aut}^0(M)$ of the group 
 of all biholomorphisms of $M$, we consider the maximal connected linear algebraic subgroup $H$ of
 $\operatorname{Aut}^0(M)$. Hence $\operatorname{Aut}^0(M)/H$ is an abelian variety.
For the identity component $Z$ of the center of a maximal compact connected subgroup $K$ of $H$, we 
take its complexification $Z_{\Bbb C}$ in $H$.
Let 
$$
\frak z := \operatorname{Lie}(Z)\quad
\text{ and  }
\quad
\frak{z}_{\Bbb C}:= \operatorname{Lie}(Z_{\Bbb C})
$$ 
be the associated Lie algebras.
Then the infinitesimal action of  $\frak{z}_{\Bbb C}$
 on $M$ lifts to an 
infinitesimal bundle action  of $\frak{z}_{\Bbb C}$  on $L$.
For
$$
V_{m} := H^0(M,\mathcal{O}(L^{\otimes m})), \qquad m = 1,2,\dots,
$$ 
we view $\frak{z}_{\Bbb C}$ as a Lie subalgebra of $\frak{sl}(V_{m})$ by taking the traceless part 
for each element of $\frak{z}_{\Bbb C}$.
In view of the infinitesimal action of $\frak{z}_{\Bbb C}$,  $V_{m}$ is expressible as 
a direct sum of $K$-invariant subspaces,
 $$
 V_{m}\; =\; \bigoplus_{\alpha =1}^{\nu_{m}} \, V_{m , \alpha },
 $$
where $ V_{m , \alpha } := \{\, \tau \in V_{m}\,;\, Y\tau = \chi^{}_{m ,\alpha} (Y) \tau 
\,\text{ for all $Y \in \frak{z}_{\Bbb C}$}\}$ 
with mutually distinct characters $\chi^{}_{m,\alpha} \in \frak{z}_{\Bbb C}^*$, 
$\alpha = 1,2,\dots, \nu_{m}$. 
Let $h$ be a $K$-invariant Hermitian metric for $L$ such that the associated first Chen form 
$\omega = c_1 (L; h)$ is K\"ahler.  Define a Hermitian metric $ \rho_{m} (h)$ for $V_{m}$ by
$$
\rho_{m} (h) (\tau, \tau' )\,:=\,\int_X\; (\tau ,\, \tau' )^{}_{h}\,\omega^{n}
\qquad \tau,\, \tau' \in V_{m},
$$
where $(\tau ,\tau' )_{h}$ denotes the pointwise Hermitian pairing of $\tau$, $\tau'$ in  terms of 
the Hermitian metric $h$. 
If $\tau$ and $\tau'$ coincide, then we write $(\tau ,\tau )_{h}$ 
simply as $|\tau |_h^{2}$.  
By setting $q:= 1/m$,
we now define 
$$
B_{m, \alpha}(h)\,:=\; n!\, q^n \,\sum_{i=1}^{n_{\alpha}}\, |\tau_{\alpha, i}|_h^2,
\qquad
B_m(h)\,:=\;\sum_{\alpha =1}^{\nu_m}\, B_{m,\alpha}(h),
$$
where  $\{\tau_{\alpha,i}\,;\, i = 1,2,\dots, n_{\alpha}\}$ is an orthonormal basis for 
the subspace $V_{m, \alpha}$ of the Hermitian vector space $(V_{m}, \rho_{m}(h))$. 
Note that, if $M$ admits an extremal K\"ahler metric $\omega_0$  in the class $c_1(L)$, then in view of \cite{C}, by choosing $K$ to be the identity component of the
group of isometries for $(M, \omega_0)$ in $H$, we may assume that the extremal K\"ahler vector field 
$\mathcal{V}$ (cf.~\cite{FM}) belongs to $\frak z$. Let $\sigma_{\omega_0}$ be the 
scalar curvature  of $\omega_0$, and define a real constant $C_0$ by
$$
C_0 \, :=\,\{2  c_1(L)^n[M]\}^{-1}\left \{ \int_M \sigma_{\omega_0}\, \omega_0^n + 
\sqrt{-1}\int_M h_0^{-1} (\mathcal{V} h_0 )\,\omega_0^n \right \}.
$$
Then by setting $\mathcal{Y}_0:= \sqrt{-1}\mathcal{V}/2$ and $\varphi_0 := 0$,  
we obtain

\medskip\noindent
{\bf Main Theorem}: {\em Suppose that $L$ admits a Hermitian metric $h_0$ such that
$\omega_0 := c_1 (L;h_0)$ is an extremal K\"ahler metric
with extremal K\"ahler vector field $\mathcal{V}\in\frak z$. Then there exist vector fields $\mathcal{Y}_k \in 
\sqrt{-1}\,\frak z$, smooth real-valued $K$-invariant functions $\varphi_k$, real constants $C_k$,  
$k = 1,2,\dots$, on $M$ such that  
$$
\,\,\sum_{\alpha =1}^{\nu_m} \{1- \chi_{m,\alpha}(\mathcal{Y}(\ell ))\} B_{m,\alpha}(h(\ell ))
\,=\, 1 + C(\ell ) + O (q^{\ell +2}),
\,\, \ell = 0,1,2,\dots,
\leqno{(1.1)}
$$
with $\mathcal{Y}_k $, $\varphi_k$, $C_k$ independent of $q$ and $\ell$, 
where
$\,\mathcal{Y}(\ell ):= \Sigma_{k=0}^{\ell}\, q^{k+2}\mathcal{Y}_k$, $h (\ell ) := h_0 \exp (- \Sigma_{k=0}^{\ell} q^k \varphi_k)$, and $C(\ell ) := \Sigma_{k=0}^{\ell}\, C_k q^{k+1}$.
}

\medskip
Here for every integer $r$, we mean by $O (q^{r})$  
a quantity whose $C^j$-norm for every nonnegative integer $j$ is bounded by 
$\kappa_j q^{r}$ for
some positive constant $\kappa_j$ independent of $q$ and $\alpha$. 
In the above Main Theorem, let $\ell \to \infty$. Then the formal expression of the left-hand side of (1.1) is called the {\it asymptotic polybalaced kernel\/} 
for $(M,L)$.

\medskip
Let $T$ 
be an arbitrary algebraic torus in $H$ satisfying $Z_{\Bbb C} \subset T$.
As a corollary of Main Theorem, 
we obtain

\medskip\noindent
{\bf Corollary}: {\em If  $M$ admits an extremal K\"ahler metric in the class $c_1 (L)$, then $(M,L)$ is asymptotically Chow-polystable relative to $T$, i.e., $(M,L^{\otimes m})$, 
$m \gg 1$, are Chow-polystable relative to $T$.
}

\medskip
In the last section, in view of \cite{MN} and a recent result of Yotsutani-Zhou \cite{YZ}, we shall discuss 
extremal K\"ahler versions of the Yau-Tian-Donaldson Conjecture from various points of view.

\medskip
For preceding related works, see Apostolov-Huang \cite{A},  Donaldson \cite{D},\cite{D1}, 
Futaki \cite{F}, Hashimoto \cite{H}, Lu \cite{L}, Ono-Sano-Yotsutani \cite{OSY}, 
Phong-Sturm \cite{PS}, Sano-Tipler \cite{ST}, Sz\'ekelyhidi \cite{S1},\cite{S2}, Tian \cite{T},\cite{T0}, 
Zelditch \cite{Ze} and  Zhang \cite{Z}. I owe much to these works.

\medskip
Parts of this paper were announced in Pacific Rim Conference on Complex  
and Simplectic Geometry  XI 
at Heifei in July, 2016.  Afterwards, during the preparation of this paper, I heard that R. Seyyedali showed 
asymptotic Chow-stability, relative to a maximal algebraic torus in $\operatorname{Aut}^0(M)$, 
for an extremal K\"ahler polarized algebraic manifold $(M,L)$.

\section{Proof of Main Theorem}

In this section, we prove Main Theorem by induction on $\ell$.
For each $K$-invariant
Hermitian metric $h$ for $L$, by setting $|\tau_{\alpha, i}|^2 := \,\tau_{\alpha, i}\,\bar{\tau}_{\alpha, i}$,
we consider
$$
\Psi_{m}(h)\; :=\; n!\, q^n \,\sum_{\alpha =1}^{\nu_m}\sum_{i=1}^{n_{\alpha}}\, |\tau_{\alpha, i}|^2,
\leqno{(2.1)}
$$
where $\{\tau_{\alpha,i}\,;\, i = 1,2,\dots, n_{\alpha}\}$ 
is an orthonormal basis for $(V_{m}, \rho_{m}(h))$ as in the introduction.
Then the left-hand side of (1.1) is the real-valued function on $X$ obtained as the
contraction 
$$
h(\ell )^m \cdot \{ (1- \mathcal{Y}(\ell ) ) \Psi_m (h (\ell ))\}
$$
of $(1- \mathcal{Y}(\ell ) ) \Psi_m (h (\ell ))$ with $h(\ell )^m$ (see \cite{MA}, (1.4.1), 
for the definition of $\mathcal{Y}(\ell )  |\tau_{\alpha, i}|^2$). Hence the proof of (1.1)
(and hence Main Theorem) is 
reduced to showing  the following for all nonnegative integers $\ell$:
$$
h(\ell )^m \cdot \{ (1- \mathcal{Y}(\ell ) ) \Psi_m (h (\ell ))\} \; =\; 1 + C(\ell ) + O (q^{\ell +2}).
\leqno{(2.2)}
$$
Note that $\mathcal{Y} (\ell )$ acts on the antiholomorphic 
section $\bar{\tau}_{\alpha, i}$ trivially.
Let $\mathcal{D}_0$ be the Lichn\'erowicz operator as defined in \cite{C}, (2.1), 
for the extremal K\"ahler manifold $(M, \omega_0)$, where we write $\omega_0$ as
$$
\omega_0 \; =\; \frac{\sqrt{-1}}{2\pi} \sum_{i,j}\, g_{i \bar{j}}\, dz^{i}\wedge d\overline{z^{j}}
$$
in terms of a system $(z^1, \dots, z^n)$ of holomorphic local coordinates on $M$.
 Let $\mathcal{S}$ be the space of all real-valued smooth $K$-invariant functions $\varphi$ on $M$ 
 such that $\int_M \varphi \omega_0^{\,n} = 0$. Since $\mathcal{D}_0$ maps $\mathcal{S}$ into itself,
the restricted operator
$$
\mathcal{D}_0 : \mathcal{S} \to \mathcal{S}
$$
is denoted also by $\mathcal{D}_0$
whose kernel in $\mathcal{S}$  
is written simply as $\operatorname{Ker}\, \mathcal{D}_0$.
Then we have an isomorphism
$$
\;\;\; e_0:\; \operatorname{Ker}\, \mathcal{D}_0\;\cong \; \frak z, \qquad \quad\varphi \leftrightarrow
e_0 (\varphi ) := \operatorname{grad}^{\Bbb C}_{\omega_0}\varphi,
\leqno{(2.3)}
$$
where $\operatorname{grad}^{\Bbb C}_{\omega_0}\varphi := (1/\sqrt{-1})\Sigma_{i,j} g^{\bar{j}i}
(\partial \varphi / \partial z^{\bar{j}}) \partial / \partial z^i$. By the inner product
$$
(\varphi ,\psi )_{\omega_0} \; :=\; \int_M \varphi \psi \,\omega_0^n,
\qquad\quad \varphi,\, \psi \in \mathcal{S},
$$
we write $\mathcal{S}$ as an orthogonal direct sum $\operatorname{Ker}\, \mathcal{D}_0 \oplus (\operatorname{Ker}\, \mathcal{D}_0)^{\perp}$. We then consider the orthogonal projection 
$\operatorname{pr}_1: \mathcal{S} \to \operatorname{Ker}\, \mathcal{D}_0$ to the first factor.
The proof of (2.2) is divided into two steps:

\medskip\noindent
{\em Step \/$1$}: In this step, we shall show that (2.2) is true for $\ell =1$. 
Note that $h(0) = h_0$. In view of Lu \cite{L}, 
the Tian-Yau-Zelditch asymptotic expansion (\cite{T}, \cite{Ze}; see also \cite{Ca}) is written in the form
$$
h_0^m \cdot  \Psi_m (h(0))\;\, (\,= \, B_m (h(0))\, )\;=\;\, 1 + \frac{\sigma_{\omega_0}}{2}q + O (q^2).
\leqno{(2.4)}
$$
By $\mathcal{Y}(0) = q^2 \sqrt{-1}\mathcal{V}/2$, we have $m\mathcal{Y}(0) = q \sqrt{-1}\mathcal{V}/2$.
Take the infinitesimal action of $\mathcal{Y}(0)$ on (2.4). Dividing it further by (2.4), we obtain
$$
q\, h_0^{-1}\sqrt{-1}\,(\mathcal{V}/2)\, h_0   \; +\;\Psi_m (h(0))^{-1}\{\mathcal{Y}(0) \Psi_m (h(0))\}
\;=\; O(q^3).
\leqno{(2.5)}
$$
On the other hand,  by \cite{M0}, p.579, we have $h_0^{-1}\sqrt{-1} (\mathcal{V}/2) h_0 = C_0 - (\sigma_{\omega_0}/2)$. Hence (2.5) is rewritten as
$$
\Psi_m (h(0))^{-1}\{\mathcal{Y}(0) \Psi_m (h(0))\}\;=\; - \,q\{C_0 - (\sigma_{\omega_0}/2)\} + O (q^3 ).
$$ 
This together with (2.4) implies that
$$
h_0^m \cdot \{\mathcal{Y}(0) \Psi_m (h(0))\}\;=\; - \,q\{C_0 - (\sigma_{\omega_0}/2)\} + O (q^2 ).
$$
Subtracting this from (2.4), we obtain the required equality:
$$
h(0)^m\cdot \{(1- \mathcal{Y}(0))\Psi_m (h(0))\} \;=\; (1+ q C_0 ) + O(q^2)\;=\; 1+ C(0) + O(q^2).
$$
{\em Step \/$2$}: Note that the left-hand side of (2.2) is
$$
\{h(\ell )^m \cdot \Psi_m (h (\ell ))\}\frac{(1- \mathcal{Y}(\ell ))\Psi_m (h(\ell ))}{\Psi_m (h(\ell ))}
\;=\; B_m (h(\ell ))\left \{ 1 - \frac{\mathcal{Y}(\ell )\Psi_m (h(\ell ))}{\Psi_m (h(\ell ))}\right\}.
$$
 For a positive integer $\ell$, by assuming the induction hypothesis 
$$
\;
B_m (h(\ell -1))\left \{ 1 - \frac{\mathcal{Y}(\ell -1)\Psi_m (h(\ell -1))}{\Psi_m (h(\ell -1))}
\right \}\, =\, 1 + C(\ell -1) + O (q^{\ell +1}),
\leqno{(2.6)}
$$
we have only to find $\mathcal{Y}_{\ell}$, $\varphi_{\ell}$ and $C_{\ell}$ such that
$\mathcal{Y} (\ell ) := \mathcal{Y} (\ell-1 ) + q^{\ell +2} \mathcal{Y}_{\ell}$, $h(\ell ) := h(\ell -1 ) e_{}^{-q^{\ell}\varphi_{\ell}}$ and $C(\ell ):= C(\ell -1 )+ C_{\ell}q^{\ell +1}$ satisfy 
$$
B_m (h(\ell ))\left \{ 1 - \frac{\mathcal{Y}(\ell )\Psi_m (h(\ell ))}{\Psi_m (h(\ell ))}
\right \}\, =\, 1 + C(\ell ) + O (q^{\ell +2}).
\leqno{(2.7)}
$$
By setting $\omega (\ell ):= c_1 (L; h(\ell ))$ and $\omega (\ell -1 ):= c_1 (L; h(\ell -1))$,
we have  $\omega (\ell ) = \omega (\ell -1 ) + (\sqrt{-1}/2\pi )q^{\ell}\partial \bar{\partial}\varphi_{\ell}$.
Let $C^{\infty}(M)^{K}_{\Bbb R}$ be the space of all real-valued $K$-invariant functions on $M$.
To each $(\mathcal{Y}_{\ell}, \varphi_{\ell}, C_{\ell})$ in 
$\sqrt{-1}\, \frak z \times C^{\infty}(M)^{K}_{\Bbb R}
\times \Bbb R$, we assign  
a real-valued $K$-invariant function
$\Phi (q; \mathcal{Y}_{\ell}, \varphi_{\ell}, C_{\ell})$ by 
\begin{align*}
&\Phi (q; \mathcal{Y}_{\ell}, \varphi_{\ell}, C_{\ell})  \, :=\, 
B_m (h(\ell  ) ) 
\left \{ 1 - \frac{\mathcal{Y} (\ell ) \Psi_m (h(\ell  ) )} {\Psi_m (h(\ell  ) )}\right \}\\
&\;\; = \; B_m (h(\ell -1 ) e_{}^{-q^{\ell}\varphi_{\ell}}) 
\left \{ 1 - \frac{(\mathcal{Y} (\ell-1 ) + q^{\ell +2} \mathcal{Y}_{\ell})\Psi_m (h(\ell -1 )e_{}^{-q^{\ell}\varphi_{\ell}}) }{\Psi_m (h(\ell -1 ) e_{}^{-q^{\ell}\varphi_{\ell}}) }\right \}.
\end{align*}
By the induction hypothesis (2.6), 
there exists a real-valued $K$-invariant function $u_{\ell}$ on $M$ such that
$$
\Phi (q; 0, 0, 0)\; \equiv\;  1 + C (\ell -1 )  + u_{\ell} q^{\ell +1},
\qquad \text{ modulo $q^{\ell +2}$}.
\leqno{(2.8)}
$$
In view of the variation formula for the scalar curvature (see for instance \cite{C}, (2.5)), we see that,
modulo $q^{\ell +2}$, 
$$
\begin{cases}
\;\;B_m (h (\ell )) - B_m (h ( \ell -1)) \; &\equiv \;\; (q/2) \{\sigma_{\omega (\ell )} - \sigma_{\omega (\ell -1)} \}\\
\;&\equiv\; \; q^{\ell +1}(- \mathcal{D}_0 + \sqrt{-1}\, \mathcal{V}) (\varphi_{\ell} /2).
\end{cases}
\leqno{(2.9)}
$$
Put $I_1 := \{\mathcal{Y}(\ell ) \Psi_m (h (\ell ))\}/\Psi_m (h (\ell ))$,
$J := \{\mathcal{Y}(\ell -1 ) \Psi_m (h (\ell ))\}/\Psi_m (h (\ell ))$
and $I_2 := \{\mathcal{Y}(\ell -1) \Psi_m (h (\ell-1))\}/\Psi_m (h (\ell -1))$. 
In view of (2.4), we have $h_0^m\cdot \Psi_m ( h (\ell )) \equiv 1$ modulo $q$.
Hence, modulo $q^{\ell +2}$,
$$
I_1 - J\;\; (\,=\; q^{\ell +2 }\{\mathcal{Y}_{\ell } \Psi_m (h (\ell ))\}/\Psi_m (h (\ell )) \,)
\;\equiv \; - \,q^{\ell +1} h_0^{-1}(\mathcal{Y}_{\ell} h_0).
\leqno{(2.10)}
$$
Note here that, by setting
$$
C \;:=\; \{c_1 (L)^n [M]\}^{-1}_{}\int_M h_0^{-1}(\mathcal{Y}_{\ell} h_0 ) \, \omega_0^n,
$$ 
we obtain $h_0^{-1} (\mathcal{Y}_{\ell} h_0) = C + e_0^{-1}(\sqrt{-1}\mathcal{Y}_{\ell} )$ 
(see for instance \cite{MA}). On the other hand, we have the following:
\begin{align*}
J - I_2\; &=\; \mathcal{Y}(\ell -1 ) \log \left \{\frac{\Psi_m (h (\ell ))}{\Psi_m (h (\ell -1 ))}\right \}\\
& =\; \mathcal{Y}(\ell -1 ) \left (q^{\ell -1}\varphi_{\ell} + \log 
\left \{\frac{h (\ell )^m\cdot\Psi_m (h (\ell ))}{h(\ell -1 )^m\cdot\Psi_m (h (\ell -1 ))}\right \}\right )\\
&=\; \mathcal{Y}(\ell -1 ) \left (q^{\ell -1}\varphi_{\ell} + \log \{ B_m (h (\ell ))/ B_m (h (\ell -1 ))\} \right ).
\end{align*}
Since $B_m (h (\ell )) \equiv B_m (h (\ell -1 )) \equiv 1$ modulo $q$, we see from (2.9) that
$\log \{ B_m (h (\ell ))/ B_m (h (\ell -1 ))\} \equiv 0$ modulo $q^{\ell +1}$.
Moreover, by $\ell \geq 1$, we obtain $\mathcal{Y}(\ell -1 )\, \equiv\, q^2 (\sqrt{-1}/2)  \mathcal{V}\,$ modulo $q^3$. It then follows that
$$
J - I_2 \; \equiv\; q^{\ell +1}\sqrt{-1}\mathcal{V}(\varphi_{\ell}/2),
\quad \text{ modulo $q^{\ell +2}$.}
\leqno{(2.11)}
$$
 On the other hand, since $h_0^m\cdot \Psi_m ( h (\ell )) \equiv 1$ modulo $q$,
we obtain 
$$
\begin{cases}
\quad  I_1 &\equiv  \;\; q^2 (\sqrt{-1}\,\mathcal{V}/2) \Psi_m (h (\ell )) /\Psi_m (h (\ell )) \\
&\equiv \;\; -\, (\sqrt{-1}/2)\,q h_0^{-1}(\mathcal{V}h_0)
\;\; \equiv\;\;  0, \qquad\text{ modulo $q$}. 
\end{cases}
\leqno{(2.12)}
$$
Now by  (2.10) and (2.11),  $I_1 - I_2 \equiv q^{\ell +1} \{\sqrt{-1}\,\mathcal{V} (\varphi_{\ell}/2) - h_0^{-1}(\mathcal{Y}_{\ell} h_0 )\} \equiv q^{\ell +1} \{\sqrt{-1}\,\mathcal{V} (\varphi_{\ell}/2) - C - e_0^{-1}(\sqrt{-1}\, \mathcal{Y}_{\ell})\}$ modulo $q^{\ell +2}$.
Thus by setting $B_1 := B_m (h (\ell ))$ and  $B_2 := B_m (h (\ell -1))$, we see from (2.9) and (2.12) the following:
$$
\begin{cases}
&{\Phi} (q; \mathcal{Y}_{\ell}, \psi_{\ell}, C_{\ell})  - 
{\Phi} (q; 0, 0, 0) \; =\; B_1 ( 1 - I_1 ) - B_2 (1 - I_2)\\
 &=\; (B_1 - B_2) (1 - I_1)  - B_2 (I_1 - I_2 ) \\
 & \equiv \; q^{\ell +1} \{(- \mathcal{D}_0 + \sqrt{-1}\, \mathcal{V}) (\varphi_{\ell} /2)\\
 &\qquad\qquad\qquad \;\;\;
  -  \sqrt{-1}\,\mathcal{V} (\varphi_{\ell}/2) + C + e_0^{-1}(\sqrt{-1}\, \mathcal{Y}_{\ell})\}\\
 & \equiv \; q^{\ell +1}\{- \mathcal{D}_0(\varphi_{\ell}/2)
+ C + e_0^{-1}(\sqrt{-1}\, \mathcal{Y}_{\ell})\}, \quad \text{ modulo $q^{\ell +2}$.}
\end{cases}\leqno{(2.13)}
$$
Since the function
$u_{\ell} -\mu_{\ell}$ belongs to $\mathcal{S}$ for $\mu_{\ell}:= \{ c_1 (L)^n [M] \}^{-1}\int_M u_{\ell} \,\omega_0^n$, 
we can write $u_{\ell}$ as a sum
$$
u_{\ell} \; =\; \mu_{\ell} + u'_{\ell} + u_{\ell}'' ,
\leqno{(2.14)}
$$
where $u'_{\ell}:= (1- \operatorname{pr}_1)(u_{\ell}- \mu_{\ell}) \in (\operatorname{Ker}\, \mathcal{D}_0)^{\perp}$ and $u_{\ell}'' := \operatorname{pr}_1(u_{\ell}- \mu_{\ell}) \in 
\operatorname{Ker}\, \mathcal{D}_0$. Let $\varphi_{\ell}$  be the unique element of $(\operatorname{Ker}\, \mathcal{D}_0)^{\perp}$ such that $\mathcal{D}_0 (\varphi_{\ell}/2) = u'_{\ell}$. Put
$$
\mathcal{Y}_{\ell} \;:=\; \sqrt{-1}\,e_0 (u_{\ell}'')
\qquad \text{ and }\qquad  C_{\ell}\;:=\; \mu_{\ell} +C.
\leqno{(2.15)}
$$
Then by (2.3), $\mathcal{Y}_{\ell} \in \sqrt{-1}\frak z$, while $C_{\ell}$ is a real constant.
Since $\mathcal{D}_0 (\varphi_{\ell}/2) = u'_{\ell}$, it follows from (2.8), (2.13), (2.14) and (2.15) that, 
modulo $q^{\ell +2}$, 
\begin{align*}
&{\Phi} (q; \mathcal{Y}_{\ell}, \psi_{\ell}, C_{\ell})\; \equiv \; {\Phi} (q; 0, 0, 0) + q^{\ell +1}\{- \mathcal{D}_0(\varphi_{\ell}/2)
+ C + e_0^{-1}(\sqrt{-1}\, \mathcal{Y}_{\ell})\}\\
&\equiv\; 1 + C(\ell -1 ) + q^{\ell +1} \{ u_{\ell} - \mathcal{D}_0(\varphi_{\ell}/2)
+ C + e_0^{-1}(\sqrt{-1}\, \mathcal{Y}_{\ell})\}\\
&= \; 1 + C(\ell -1 ) + q^{\ell +1}\{ (u'_{\ell} - \mathcal{D}_0(\varphi_{\ell}/2)) + (\mu_{\ell} + C )
+ (u_{\ell}'' +  e_0^{-1}(\sqrt{-1}\, \mathcal{Y}_{\ell}))\}\\
&= \;  1 + C(\ell -1 ) + q^{\ell +1} C_{\ell} \; \equiv \; 1 + C(\ell ),
\end{align*}
which shows (2.7), as required.
\qed

\medskip\noindent
{\em Remark\/} 2.16: The preceding work in \cite{M0}, Theorem B, is obtained from Main Theorem above
by replacing the left-hand side of (1.1) by
$$
\,\,\sum_{\alpha =1}^{\nu_m} \exp \{- \chi_{m,\alpha}(\mathcal{Y}(\ell ))\}\, B_{m,\alpha}(h(\ell )).
\leqno{(2.17)}
$$ 
The point is the following:
The coefficients 
$1- \chi_{m,\alpha} (\mathcal{Y}(\ell ))$, $\alpha = 1,2,\dots, \nu_m$, in (1.1) are linear 
in $\mathcal{Y} (\ell )$, while the coefficients 
$\exp \{- \chi_{m,\alpha}(\mathcal{Y}(\ell ))\}$, $\alpha = 1,2,\dots, \nu_m$, in (2.17) aren't.
This linearity is essential in the proof of the asymptotic relative Chow-polystability.

\section{Relative Chow-polystability}

In this section, we fix an algebraic torus $T$ in $\operatorname{Aut}^0(M)$.
Then for $K$ in the introduction, replacing $T$ by its conjugate group, 
we may assume that the maximal compact subgroup $T_c$ of $T$ sits in $K$.
Put $\frak t_c:=\operatorname{Lie}(T_c)$.
Note that the infinitesimal action of the Lie algebra 
$\frak{t}$ 
of $T$ lifts to an 
infinitesimal bundle action of $\frak{t}$ on $L$.
For each positive integer $m$, let $V_{m} := H^0(M,\mathcal{O}(L^{\otimes m}))$,
and we view $\frak t$ as a Lie subalgebra of $\frak{sl}(V_{m})$ by considering the traceless part.
Define a symmetric bilinear form $\langle \,\, ,\, \rangle_{m}$ on $\frak{sl}(V_{m})$ by
$$
\langle X, Y \rangle_m^{} \;:=\; \operatorname{Tr}(XY)/m^{n+2},
\qquad X,Y\in \frak{sl}(V_{m}).
$$
Let $\frak{z}_{m}$ be the centralizer of $\frak t$ in $\frak{sl}(V_{m})$, 
and $\frak{g}_{m}$ be the orthogonal complement of $\frak t$ in $\frak{z}_{m}$,
i.e.,
$$
\begin{cases}
\;\;\frak{z}_{m} & := \{ X \in \frak{sl}(V_{m})\,;\, [X,Y] = 0 \;\text{ for all $Y \in \frak{t}$} \},\\
\;\;\frak{g}_{m} & := \{ X \in \frak{z}_{m}\,;\, \langle X, Y\rangle_{m} = 0 \;\text{ for all $Y \in \frak{t}$} \}.
\end{cases}
\leqno{(3.1)}
$$
Let $Z_{m}$ and $G_{m}$ be the connected reductive algebraic subgroups in
$\operatorname{SL}(V_{m})$
associated to $\frak{z}_{m}$ and $\frak{g}_{m}$, respectively.
 By the infinitesimal $\frak{t}$-action  on $V_{m}$, we can write $V_{m}$ as a direct sum 
 of $T_c$-invariant subspaces,
 $$
 V_{m}\; =\; \bigoplus_{\gamma =1}^{\eta_{m}} \, \hat{V}_{m , \gamma },
 \leqno{(3.2)}
 $$
where $ \hat{V}_{m , \gamma } = \{ \sigma \in V_{m}\,;\, Y\sigma = \hat{\chi}^{}_{m ,\gamma} (Y) \sigma 
\;\text{ for all $Y \in \frak t $}\}$ with mutually distinct characters $\hat{\chi}^{}_{m ,\gamma} \in \frak{t}^*$, 
$\gamma = 1,2,\dots, \eta_{m}$. We now consider the algebraic subgroup $R_{m}$ of $\operatorname{SL}(V_{m})$ defined by
$$
R_{m} \;:=\; \prod_{\gamma =1}^{\eta_{m}}\,
\operatorname{SL}(\hat{V}_{m , \gamma }),
$$
where each $\operatorname{SL}(\hat{V}_{m , \gamma })$ fixes $\hat{V}_{m , \gamma' }$ if 
$\gamma' \neq \gamma$. Then the centralizer $H_{m}$ of $R_{m}$ in $\operatorname{SL}(V_{m})$ 
consists of all diagonal matrices in $\operatorname{SL}(V_{m})$ acting on each $\hat{V}_{m , \gamma }$, 
$\gamma = 1,2,\dots, \eta_m$, by 
constant scalar multiplication. Note that $\frak{t}$ viewed as a Lie subalgebra of $\frak{sl}(V_{m})$
 sits in the Lie algebra $\frak{h}_{m}$ of $H_{m}$.  Let
$$
\frak{t}_{m}^{\perp}\; := 
\{ X \in \frak{h}_{m}\,;\, \langle X, Y\rangle_{m} = 0 \;\text{ for all $Y \in \frak{t}$} \}.
$$ 
be the orthogonal complement of $\frak{t}$ in $\frak h_{m}$.
Let $T_{m}^{\perp}$ denote the corresponding algebraic torus sitting in $H_{m}$.
Since $Z_{m} = H_{m}\cdot R_{m}$, it follows that
$$
G_{m} \; =\; T_{m}^{\perp} \cdot R_{m}.
$$
Let $M_{m}$ be the image of $M$ under the Kodaira embedding $\Phi_{m}: M \hookrightarrow 
\Bbb P^* (V_{m} )$ associated to the complete linear system $|L^{\otimes m}|$ on $M$.
For the degree $d_{m}$ of $M_{m}$ in $\Bbb P^*(V_{m})$, we consider the space
$$
W_{m} \;:=\;\{\operatorname{Sym}^{d_{m}}(V_{m})\}_{}^{\otimes n+1},
$$
where $\operatorname{Sym}^{d_{m}}(V_{m})$ denotes the $d_{m}$-th symmetric tensor 
product of $V_{m}$. For the dual space $W_{m}^*$ of $W_{m}$, let 
$0 \neq \tilde{M}_{m} \in W_{m}^*$ denote the Chow form for the irreducible reduced 
algebraic cycle $M_{m}$ on $\Bbb P^*(V_{m} )$, so that the associated 
point $[ \tilde{M}_{m} ]$ in $\Bbb P^*(W_{m} )$ is the Chow point for $M_{m}$.

\medskip
For relative stability, we showed in \cite{M1} (see also \cite{M3})
that, if $c_1(L)$ admits an extremal K\"ahler metric, then
the orbit $R_m \cdot \tilde{M}_m$ is 
closed in $W_{m}^*$ for all sufficiently large $m$.
In \cite{S1} (see also \cite{S2}), Sz\'ekelyhidi introduced the following stronger stability concept:

\medskip\noindent
{\em Definition \/$3.2$.} (1) A polarized algebraic manifold $(M,L^{\otimes m})$ is called {\it Chow-polystable relative to $T$}, if 
$G_{m}\cdot \tilde{M}_{m}$ is closed in $W_{m}^*$.

\smallskip\noindent
(2) $(M,L)$ is called {\it asymptotically Chow-polystable relative to $T$}, if $(M,L^{\otimes m})$ is 
Chow-polystable relative to $T$ for $m \gg 1$.

\medskip\noindent
{\em Definition \/$3.3$.} A polarized algebraic manifold $(M,L^{\otimes m})$ is called {\it Chow-stable relative to $T$}, if the following conditions are satisfied:

\smallskip\noindent
(1) $G_{m}\cdot \tilde{M}_{m}$ is closed in $W_{m}^*$.
\newline
(2) The isotropy subgroup of $G_{m}$ at $[\tilde{M}_{m}]$ is finite.

\section{Proof of Corollary}

In this section, using the same notation as in the preceding sections, 
we assume that $M$ admits an extremal K\"ahler metric $\omega_0 = c_1 (L;h_0)$ in the class $c_1 (L)$,
where $h_0$ is a Hermitian metric for $L$.
Following the arguments in \cite{M2}, we shall show that $(M, L^{\otimes m})$ are Chow-polystable 
relative to $T$ for $m \gg 1$.
As in the introduction, we may assume that $K$ is the identity component of the group of isometries for $(M,\omega_0)$.
Put $\frak k := \operatorname{Lie}(K)$.
Let $[n/2]$ be the largest integer which does not exceed $n/2$. By applying Main Theorem  to
$$
\ell \; :=\; [n/2] +3,
$$
we obtain a $K$-invariant K\"ahler metric $\omega (\ell ) = c_1 (L; h(\ell ))$ in the class $c_1 (L)$
such that (1.1) holds.
 For the compact group
 $$
K_{m} := \operatorname{SU}(V_{m};\rho_{m}(h (\ell ))) \cap G_{m},
 $$
we  can view $G_{m}$ as its complexification. 
Then for the $G_{m}$-action on $W^*_{m}$, 
the isotropy subgroup 
of $K_{m}$ at $\tilde{M}_{m}$ has the Lie algebra $\frak k_0$ sitting in $\frak k$.
Since $Z \subset T$, the isotropy subgroup of $K_{m}$ at $[\tilde{M}_{m}]$
has the same Lie algebra $\frak k_0$.
For $\frak g_{m} : = \operatorname{Lie}(G_{m})$, 
we define the Lie subalgebras $\frak p_{m}$ and $\frak p$ by
$$
\frak p_{m} := \sqrt{-1}\, \frak k_{m}
\quad \text{ and } \quad 
\frak p := \sqrt{-1}\, \frak k_0,
$$
where $\frak k_{m}:= \operatorname{Lie}(K_{m})$.
Put $n_{\gamma} :=\dim \hat{V}_{m , \gamma }$.
By choosing an orthonormal basis $\{ \,\sigma_{\gamma, i}\,;\, i=1, 2, \dots, n_{\gamma}\,\}$ 
for $(\hat{V}_{m,\gamma}, \rho_{m}(h (\ell )))$, we set
$$
j(\gamma , i ) \,:= \, i\, +\, \sum_{\beta =1}^{\gamma -1}\, n^{}_{\beta},
\qquad i = 1,2,\dots, n^{}_{\gamma};\; \gamma = 1,2,\dots, \eta_{m},
\leqno{(4.1)}
$$
where the right-hand side denotes $i$ in the special case $\gamma =1$.
Let $m \gg 1$.
Then for each $\gamma$ and $i$ as above,
we put
$$
\hat{\sigma}_{\gamma , i} \; 
:=\,\sqrt{1 - \hat{\chi}_{m,\gamma}(\mathcal{Y}(\ell ))}\,{\sigma}_{\gamma , i}.
\leqno{(4.2)}
$$
By writing $\hat{\sigma}_{\alpha , i}$, ${\sigma}_{\alpha , i}$
as $\hat{\sigma}_{j(\alpha , i )}$, ${\sigma}_{j(\alpha , i )}$, by abuse of terminology,
we have bases 
$$
\{\hat{\sigma}_1, \hat{\sigma}_2, \dots, \hat{\sigma}_{N_{m}}\}\quad \text{ and }\quad
\{{\sigma}_1, {\sigma}_2, \dots, {\sigma}_{N_{m}}\},
\leqno{(4.3)}
$$
respectively,
for $(V_{m}, \rho_{m}(h (\ell )))$.
Let $\Phi_{m} : M \hookrightarrow \Bbb P^{N_m -1}(\Bbb C )\;( = \Bbb P^*(V_{m}))$ be 
the associated Kodaira embedding
defined by
$$
\Phi_{m} ( p ) \,:=\, (\hat{\sigma}_1 ( p ) : \hat{\sigma}_2 ( p ): \dots : \hat{\sigma}_{N_{m}} ( p ) ),
\qquad p \in M.
\leqno{(4.4)}
$$
Here $V_m$ and $\Bbb C^{N_m}$ are identified 
by the basis $\{\hat{\sigma}_1, \hat{\sigma}_2, \dots, \hat{\sigma}_{N_{m}}\}$. 
Let $g_{\operatorname{euc}}$ be the Euclidean metric for 
the space
$\Bbb C^{N_m}= \{ (z_1, z_2, \dots, z_{N_m}) \}$.
Define the Fubini-Study form $\omega_{\operatorname{FS}}$ on 
$\Bbb P^*(V_{m}) \; (=\{ (z_1:z_2: \dots :z_{N_m}\})$
by
$$
\omega_{\operatorname{FS}}\,:=\, (\sqrt{-1}/2\pi ) 
\partial\bar{\partial} \log (\Sigma_{j = 1}^{N_{m}} |z_{j} |^2 ).
$$
For each $X \in \frak p_{m}$, let $\mathcal{V}_X$ be the associated holomorphic vector field on $\Bbb P^*(V_{m})$. We then have a unique real-valued function $\varphi_X$ on  $\Bbb P^*(V_{m})$ 
satisfying 
$$
\int_{\Bbb P^{N_m -1}(\Bbb C )} \varphi_X \,\omega_{\operatorname{FS}}^{N_{m} -1} = 0
\quad \text{ and }\quad
i_{\mathcal{V}_X} (\omega_{\operatorname{FS}}/m ) \; =\; (\sqrt{-1}/2\pi ) \bar{\partial} \varphi_X.
$$
Let us consider the  real-valued function $\zeta = \zeta (x)$ on $\Bbb R$ defined by
$\zeta (x) := x (e^x + e^{-x} )/ (e^x - e^{-x} ), \; x\in \Bbb R$.
In view of $M_{m} = \Phi_{m}(M)$,
we define a positive semidefinite $K$-invariant inner product $(\;,\, )_{m}$ on $\frak p_{m}$ by
$$
(X, Y)_{m} \;=\; \sqrt{-1} \int_{M_{m}} \partial\varphi^{}_Y \wedge \bar{\partial}\varphi^{}_X\wedge n \omega_{\operatorname{FS}}^{n-1},
\qquad X,\, Y \in \frak p_{m}.
$$
Then this inner product is positive definite when restricted to $\frak p$.
Hence as a vector space, $\frak p_{m}$ is written as an orthogonal direct sum 
$\frak p \oplus \frak p^{\perp}$, where $\frak p^{\perp}$ is the orthogonal complement of $\frak p$ 
in $\frak p_{m}$. Define an open neighborhood 
$$
U_{m} \,: = \, \{ \, X \in \frak p^{\perp}\,;\, \zeta (\operatorname{ad} X) \frak p \cap \frak p^{\perp} = \{0\}\,\}
$$
of the origin in $\frak p^{\perp}$. Let $0 \neq X \in\frak p^{\perp}$.  Since $X$ belongs to $\frak g_m$, 
by choosing a suitable orthonormal basis $\{ \,\sigma_{\gamma, i}\,;\, i=1, 2, \dots, n_{\gamma}\,\}$ 
for $(\hat{V}_{m,\gamma}, \rho_{m}(h (\ell )))$, we obtain real constants $b_j$ such that
$$
X \hat{\sigma}_j \;=\; b_j \hat{\sigma}_j,
\qquad j= 1,2,\dots ,N_m,
\leqno{(4.5)}
$$
where $\{\hat{\sigma}_1, \hat{\sigma}_2, \dots, \hat{\sigma}_{N_{m}}\}$ is the basis
for $(V_{m}, \rho_{m}(h (\ell )))$ as in (4.3). 
Define a real one-parameter subgroup $\lambda_X : \Bbb R_+ \to G_m$ of $G_m$ by
$\lambda_X (e^t ) \, :=\, \exp (tX)$, $ t \in \Bbb R$.
We then consider a real-valued function $f_{X,m} (t) $ on $\Bbb R$ defined by
$$
f_{X,m} (t) \,:=\, \log \| \lambda_X ( e^t ) \cdot \tilde{M}_{m} \|_{\operatorname{CH}(g_{\operatorname{euc}})},
\qquad t \in \Bbb R,
\leqno{(4.6)}
$$
where $W _{m}^*\owns w \mapsto \|w\|_{\operatorname{CH}(g_{\operatorname{euc}})} \in \Bbb R_{\geq 0}$ 
is the Chow norm by Zhang \cite{Z} (see also \cite{M0}). Put $\dot{f}_{X,m}(t) := (d/dt) {f}_{X,m}$ and 
$\ddot{f}_{X,m}(t) := (d^2/dt^2){f}_{X,m}$.
Let 
$$
\delta_0 := \frac{q}{ \sqrt{\Sigma_{j=1}^{N_m} b_j^{\,2}}}.
$$
Then by Lemma 3.4 in \cite{M2}, the proof of Corollary is reduced to showing that,
if $m \gg 1$, then for every $0\neq X \in \frak p^{\perp}$, 
$$
\dot{f}^{}_{X,m}(t_{m}) \; =\; 0 \; < \; \ddot{f}_{X,m}(t_{m})
\qquad \text{and} \qquad  t_{m} \cdot X \in U_{m},
\leqno{(4.7)}
$$
where $t_{m}$ is a suitable real number satisfying $|t_{m}| < \delta_0$.
Now the proof is divided into the following three steps:

\medskip\noindent
{\em Step \/$1$}. Put $\bar{b} := \max_j |b_j|$. Let $t$ be an arbitrary real number 
satisfying 
$$
|t| \;  \leq \;\delta_0.
$$
It then follows that $|t| \leq q/\bar{b}$. Put $\lambda_t := \lambda_X (e^t)$ 
and $M_{m,t}:= \lambda_t ( M_m)$.
Let $T\Bbb P^* (V_m)$ and $TM_{m,t}$ denote the holomorphic tangent bundles of
$\Bbb P^* (V_m)$ and  $M_{m,t}$, respectively.
From now on, by $C_i$, $ i =0,1,2,\dots$, we mean
positive real constants independent of 
 the choice of the triple $(m,t, X)$.  
Let $m \gg 1$. Then for each integer $k \geq 0$, the argument in \cite{M2}, Step 1, 
shows that
$$
\| \omega_0 - (1/m) \Phi_m^* \lambda_t^* \omega_{\operatorname{FS}} \|_{C^k (M, \omega_0 )}
\; \leq \; C_0.
\leqno{(4.8)}
$$
Here $C_1$ possibly depends on $k$.
 From now on, $X$ viewed as a holomorphic vector field on $\Bbb P^* (V_m)$ 
will be denoted by $\mathcal{X}$.
Metrically, we identify 
the normal bundle of $M_{m,t}$ in $\Bbb P^*(V_m)$
with the subbundle $TM_{m,t}^{\perp}$ of $T\Bbb P^* (V_m)_{|M_{m,t}}$
obtained as the the orthogonal complement of $TM_{m,t}$
in  $T\Bbb P^* (V_m)_{|M_{m,t}}$.
Hence $T\Bbb P^* (V_m)_{|M_{m,t}}$ is differentiably written as the direct sum
$TM_{m,t}\oplus 
TM_{m,t}^{\perp}$.
Associated to this, the restriction $\mathcal{X}_{|M_{m,t}}$ of $\mathcal{X}$
to $M_{m,t}$ is written as
$$
\mathcal{X}_{|M_{m,t}} \; =\; \mathcal{X}_{TM_{m,t}}
\oplus \mathcal{X}_{ TM_{m,t}^{\perp} },
$$
where $\mathcal{X}_{TM_{m,t}}$ and $\mathcal{X}_{ TM_{m,t}^{\perp} }$ are $C^{\infty}$ 
sections of $TM_{m,t}$ and $TM_{m,t}^{\perp}$, respectively. Then the second derivative 
$\ddot{f}_{X,m}(t)$ is given by
$$
 \ddot{f}_{X,m} (t) \; =\; \int_{M_{m,t}} 
 | \mathcal{X}_{ TM_{m,t}^{\perp} } |^2_{\omega_{\operatorname{FS}}}\,
 \omega_{\operatorname{FS}}^n
 \; \geq \; 0.
$$
In view of (4.8), it follows from the argument of Phong and Sturm \cite{PS} that (cf.~\cite{M2}, p.235; see also \cite{D})
$$
 \begin{cases}
&\int_{M_{m,t}}  | \mathcal{X}_{TM_{m,t}^{\perp}}  |^2_{\omega_{\operatorname{FS}}}\,
 \omega_{\operatorname{FS}}^n \; \geq \; C_1 q \int_{M_{m,t}}  | \mathcal{X}_{TM_{m,t}} |_{\omega_{\operatorname{FS}}}^2\,
 \omega_{\operatorname{FS}}^n,\\
&\ddot{f}_{X,m} (t) \; \geq \; C_2\, q \int_{M_{m}} | \mathcal{X}_{|M_{m}}|^2_{\omega_{\operatorname{FS}}}\,\omega_{\operatorname{FS}}^n\;
\geq  \; C_2\, q \int_M \Theta \,\Phi_m^*\omega_{\operatorname{FS}}^n, 
\end{cases} \leqno{(4.9)}
 $$
 where $\Theta := (\Sigma_{j=1}^{N_m}\, |\hat{\sigma}_j |^2 )^{-2}\{ (\Sigma_{j=1}^{N_m}\, |\hat{\sigma}_j |^2 ) (\Sigma_{j=1}^{N_m}\, b_j^{\,2} |\hat{\sigma}_j |^2) - (\Sigma_{j=1}^{N_m}\, b_j |\hat{\sigma}_j |^2)^2\} \, \geq \,0$. Moreover, by (3.4.2) in \cite{Z} (see also \cite{M0}),
$$
\dot{f}_{X,m} (0) \; =\; (n+1)\int_M \frac{\Sigma_{j=1}^{N_m}\, b_j |\hat{\sigma}_j |^2_{h(\ell )}}{\Sigma_{j=1}^{N_m}\, |\hat{\sigma}_j |_{h(\ell )}^2}\,\Phi_m^*\omega_{\operatorname{FS}}^n.
\leqno{(4.10)}
$$
 
 \medskip\noindent
{\em Step \/$2$}.  For
the basis $\{\hat{\sigma}_1, \hat{\sigma}_2, \dots, \hat{\sigma}_{N_m}\}$ 
as above satisfying (4.2), we can write $\Psi_m (h (\ell ) )$ in (2.1) applied to $h = h(\ell )$ in the 
following form:
$$
\Psi_m ( h(\ell )) \; =\; n! \, q^n \sum_{j=1}^{N_m} \, |{\sigma}_j |^2
\;=\;  n! \, q^n \sum_{\gamma =1}^{\eta_m} 
\sum_{i=1}^{n_{\gamma}} |{\sigma}_{\gamma ,i} |^2.
$$
Then the left-hand side of (1.1) is expressible as
\begin{align*}
h (\ell )^m \cdot \{(1 - \mathcal{Y}(\ell ) ) \Psi_m (h (\ell )) \}
\; &=\; n! \, q^n \sum_{\gamma =1}^{\eta_m} 
\sum_{i=1}^{n_{\gamma}}\,\{1 - \hat{\chi}_{m,\gamma}(\mathcal{Y}(\ell ))\}\,|{\sigma}_{\gamma , i}|^2_{h(\ell )} \\
&=\; n! \, q^n \sum_{j =1}^{N_m} 
\, |\hat{\sigma}_{j}|^2_{h(\ell )},
\end{align*}
and hence (1.1) is written as
$$
n! \, q^n\sum_{j =1}^{N_m} 
\, |\hat{\sigma}_{j}|^2_{h(\ell )} \; =\;  1 + C (\ell ) + O (q^{\ell + 2} ).
\leqno{(4.11)}
$$
By taking $(\sqrt{-1}/2\pi){\partial}\bar{\partial} \log$ of both sides of (4.11), we obtain
$$
\omega_{\operatorname{FS}} \, -\, m \omega (\ell ) \; =\; O(q^{\ell +2} ).
\leqno{(4.12)}
$$
For each $\gamma$ and $i$, we put $a_{\gamma, i}:= 
\hat{\chi}_{m,\gamma} (\mathcal{Y}(\ell ) )$, where $a_{\gamma, i}$ is obviously independent of 
the choice of $i$.
Then by $\hat{\chi}_{m,\gamma} = O (q^{-1})$ and $\mathcal{Y}(\ell ) = O (q^2)$, we see that
$a_{\gamma, i} = O(q)$, i.e., $|a_{\gamma, i}| \leq \kappa q$ for some 
positive constant $\kappa$ independent of the choice of $(m, \gamma, i)$.
We also obtain
$$
|\hat{\sigma}_{\gamma , i}|_{h(\ell )}^2 \;
= \;\{1 - \hat{\chi}_{m,\gamma}(\mathcal{Y}(\ell ))\}\,|{\sigma}_{\gamma , i}|_{h(\ell )}^2
\;=\; (1- a_{\gamma, i} ) \,|{\sigma}_{\gamma , i}|_{h(\ell )}^2.
\leqno{(4.13)}
$$
In terms of (4.1), we write $a_{\gamma, i}$ 
as $a_{j(\gamma, i )}$.
Note that $X$, as an element of $\frak p^{\perp}$,  belongs to $\frak g_m$, so that $X$ sits in $\frak{sl}(V_m)$.
This together with (3.1) implies
$$
\sum_{j=1}^{N_m} \, b_j \; =\;\sum_{j=1}^{N_m} \,b_j a_j \; =\; 0. 
$$
Hence, by setting 
$$
I_m\, :=\, \Sigma_{j=1}^{N_m}\, b_j (1- a_j )  |\sigma_j |^2_{h(\ell )}
\; (=\Sigma_{j=1}^{N_m}\,  b_j |\hat{\sigma}_j |^2_{h(\ell )})
\leqno{(4.14)}
$$ 
and $\xi_1 := 
(n+1)!\{1+ C (\ell )\}^{-1}$,  we see from (4.10), (4.11), (4.12) and (4.13) that $\dot{f}_{X,m}(0)$ 
is expressible as
$$
\;
\begin{cases}
& (n+1) \int_M  \{n! q^n\Sigma_{j=1}^{N_m}\,   |\hat{\sigma}_j |^2_{h(\ell )}\}^{-1} (n! q^n I_m)\, \Phi_m^*\omega_{\operatorname{FS}}^n\\
&=\; (n+1)!\int_M \{1+ C (\ell ) + O (q^{\ell +2} )\}^{-1}\, I_m \, \{ \omega (\ell ) + O(q^{\ell +3})\}^n\\
&= \; \xi_1\int_M \{1+ O(q^{\ell +2})\}\,I_m \,
\omega (\ell ) ^n\\
&=\; \xi_1
 \{   ( \Sigma_{j=1}^{N_m}\, b_j  -\Sigma_{j=1}^{N_m}\, b_j a_j   ) +  
\int_M  O(q^{\ell +2})\,I_m \,
\omega (\ell ) ^n
 \}\\
&=\; \int_M  O(q^{\ell +2})\,I_m \,
\omega (\ell ) ^n.
\end{cases}
\leqno{(4.15)}
$$
Then by (4.12) and  (4.15) together with the second line of (4.9), we see that
$$
\;\;\;\;\;\begin{cases}
\;\;\dot{f}_{X,m} (\delta_0 ) &\geq \; \dot{f}_{X,m}(0)\, +\, C_2\, \delta_0\, q \int_M  \Theta \,\Phi_m^*\omega_{\operatorname{FS}}^n\\
&\geq \; \int_M \{O(q^{\ell +2}) \,I_m   \,+\,  C_3\, \delta_0\, q^{1-n}\, \Theta \} \,\omega (\ell )^n,\\
\;\;\dot{f}_{X,m} (-\delta_0 ) &  \leq  \;  \dot{f}_{X,m}(0)\, -\, C_2\, \delta_0\, q \int_M  \Theta \,\Phi_m^*\omega_{\operatorname{FS}}^n            \\
&\leq\; \int_M \{O(q^{\ell +2}) \,I_m   \,-\,  C_3\, \delta_0\, q^{1-n}\, \Theta \} \,\omega (\ell )^n.
\end{cases}
$$
Now as in  \cite{M2}, Remark 4.25, the inclusion $t_m X \in U_m$ follows from 
$|t_m| < \delta_0$, where in the proof, we use the basis $\{\hat{\sigma}_1, \dots, \hat{\sigma}_{N_m}\}$
in place of $\{{\sigma}_1, \dots, {\sigma}_{N_m}\}$.
Hence, in order to prove (4.7), it suffices to show the following 
for $m \gg 1$:
$$
\int_M \, \Theta \,\omega (\ell )^n \; > \; 0  \quad \text{and}\quad
 \frac{ \int_M \,q^{\ell+2}\, |I_m|  \, \omega (\ell )^n }{\int_M \,\delta_0\, q^{1-n}\,\Theta \,\omega (\ell )^n}
\;\ll\; 1.
\leqno{(4.16)}
$$

\smallskip\noindent
{\em Step \/$3$}.  Put $B_0 := \Sigma_{j=1}^{N_m}\, |\hat{\sigma}_j|_{h (\ell )}^2$,
$B_1 := \Sigma_{j=1}^{N_m}\, b_j^{\,2}|\hat{\sigma}_j|_{h (\ell )}^2$,
$B_2 := \Sigma_{j=1}^{N_m}\, b_j|\hat{\sigma}_j|_{h (\ell )}^2$.
By setting $\theta_1 := \int_M (B_1/B_0)\,\omega (\ell )^n$ and 
$\theta_2 := \int_M (B_2/B_0)^2\,\omega (\ell )^n$, we obtain
$$
\int_M \, \Theta \,\omega (\ell )^n\; =\; \theta_1 - \theta_2,
$$
where both $\theta_1$ and $\theta_2 $ are obviously nonnegative. Moreover, for $m \gg 1$, 
$$
\begin{cases}
\;\;\theta_1 &=\;\;\int_M \,{(n!q^n B_0)}^{-1}{n!q^n B_1}\, \omega (\ell )^n\\
&=\;\; \int_M \,\{1+ C(\ell )+ O(q^{\ell +2})\}^{-1}{n!q^n B_1}\, \omega (\ell )^n\\
&\geq \; \;n! q^n (1- \varepsilon )\Sigma_{j=1}^{N_m}\,b_j^{\,2},
\end{cases}
\leqno{(4.17)}
$$
where $\varepsilon > 0$ is a small  real constant independent of  $m$.
On the other hand, we see from (4.14) that 
$$
\int_M \,q^{\ell+2}\, |I_m|  \, \omega (\ell )^n\; \leq \; q^{\ell+2}\,
 \Sigma_{j=1}^{N_m}\,  | b_j (1-a_j) |.
 \leqno{(4.18)}
$$
Now the following cases are possible:

\medskip\noindent
Case 1:\, $\theta_1 > 2 \theta_2$,  \qquad Case 2:\, $\theta_1 \leq  2 \theta_2$.

\medskip
Suppose Case 1 occurs. Then $\theta_1 - \theta_2  = (1/2) \theta_1 + (1/2) ( \theta_1 - 2 \theta_2) >0$,
i.e., the first inequality in (4.16) holds. Let $m \gg 1$. Then by $a_j = O(q)$, 
$$
0 \;<\; 1 - a_j \; <\; 1 + \varepsilon,
$$
where $\varepsilon$ is as above.
Let L.H.S. denotes the left-hand side of the second inequality in (4.16).
Then in view of
(4.17) and (4.18), we obtain 
\begin{align*}
&\text{L.H.S.} \; \left ( = \;\frac{ \int_M \,q^{\ell+2}\, |I_m|  \, \omega (\ell )^n }{\int_M \,\delta_0\, q^{1-n}\,\Theta \,\omega (\ell )^n}\,\right )
\; \leq \; \frac{q^{\ell+2}\, \Sigma_{j=1}^{N_m}\,| b_j (1-a_j) |}{\delta_0\,q^{1-n}\, (\theta_1 - \theta_2 )}
 \\
& \leq \; \frac{q^{\ell+2}\, \Sigma_{j=1}^{N_m}\, |b_j (1-a_j) |}{(1/2)\,\delta_0\,q^{1-n}\, \theta_1}
\leq \; \frac{q^{\ell +2}(1+\varepsilon )\Sigma_{j=1}^{N_m}|b_j |}{(1/2)\,\delta_0\,q\,n!  (1- \varepsilon )\Sigma_{j=1}^{N_m}\,b_j^{\,2}}\\
&= \;\frac{q^{\ell }(1+\varepsilon )\Sigma_{j=1}^{N_m}|b_j |}{(1/2)\,n!  (1- \varepsilon )\sqrt{\Sigma_{j=1}^{N_m}\,b_j^{\,2}}}\; \leq \; \frac{q^{\ell }(1+\varepsilon ) N_m^{1/2}}{(1/2)\,n!  (1- \varepsilon )},
\end{align*}
where the Schwarz inequality 
$\Sigma_{j=1}^{N_m} |b_j | \leq  N_m^{1/2} \sqrt{\Sigma_{j=1}^{N_m} b_j^{\,2}}$
is used in the last inequality. Hence by $N_m = O (m^n)$, it follows that
$$
\text{L.H.S.} \; \leq \; O(q^{\ell - \frac{n}{2}}).
$$
In view of the definition $\ell := [n/2] + 3$, we have $\ell - \frac{n}{2} >0$, and therefore 
$\text{L.H.S.} \ll 1$, as required. Thus in this case (4.7) holds.

\medskip
Next we consider the situation where Case 2 occurs. Note that, for $0 \neq X \in \frak p^{\perp}$, 
we can write 
$$
\Phi_m^*\varphi_X\; =\; \frac{\Sigma_{j=1}^{N_m} b_j |\hat{\sigma}_j|^2}{m\Sigma_{j=1}^{N_m}  |\hat{\sigma}_j|^2}\; =\; \frac{B_2}{mB_0}.
\leqno{(4.19)}
$$
Let $c_X$ be the real constant such that the function $\phi_X := c_X + \Phi_m^*\varphi_X$ on $M$ satisfies
$\int_M \phi_X \, \tilde{\omega}^n  =0$, where $\tilde{\omega} :=\Phi_m^*(q \omega_{\operatorname{FS}})$.  Then
$$
\|\phi_X\|^2_{L^2(M,\, \tilde{\omega})}
\;\leq \; C_4\|\bar{\partial}\phi_X\|^2_{L^2(M, \, \tilde{\omega})}
\;=\; C_4 \|\Phi_m^* \mathcal{X}_{TM_m}\|^2_{L^2(M, \, \tilde{\omega})}
$$
for some $C_4$, while by the first inequality in (4.9) applied to $t=0$, 
$$
 \|\Phi_m^* \mathcal{X}_{TM_m}\|^2_{L^2(M, \, \tilde{\omega})} \; \leq \; C_1^{-1}q^{-1}
 \|\Phi_m^* \mathcal{X}_{TM^{\perp}_m}\|^2_{L^2(M, \, \tilde{\omega})}.
$$
In view of these inequalities, we obtain
$$
\|\phi_X\|^2_{L^2(M,\, \tilde{\omega})} \; \leq \; C^{}_4\,C_1^{-1}q^{-1}
 \|\Phi_m^* \mathcal{X}_{TM^{\perp}_m}\|^2_{L^2(M, \, \tilde{\omega})}.
 \leqno{(4.20)}
$$
Let $ m \gg 1$.   By (4.10) and (4.19),
$$
|c_X| \;= \left (\int_M \tilde{\omega}^n\right )^{-1} 
\left |\int_M (\Phi_m^*\varphi_X )\, \tilde{\omega}^n  \right |
\;= \;
\frac{q^{n+1}| \dot{f}_{X,m}(0)  |}{(n+1)c_1 (L)^n[M]}. 
\leqno{(4.21)}
$$
In view of (4.12) and (4.15), there exist $C_5$ and $C_6$ satisfying
$$
| \dot{f}_{X,m}(0)  |\,\leq\, C_5\,q^{\ell +2}\| I_m \|_{L^1(M, \omega (\ell ))}
\, \leq \, C_6\,q^{\ell +2}\| I_m \|_{L^2(M, \tilde{\omega} )},
$$
while by (4.11), (4.14) and (4.19), we obtain $C_7$  
such that
$$
\| I_m \|_{L^2(M, \tilde{\omega} )}\; 
\leq \; C_7\,m^{n+1} \| \Phi_m^*\varphi_X\|_{L^2(M, \tilde{\omega} )}.
$$
Hence for some $C_8$, it follows that
$$
| \dot{f}_{X,m}(0)  | \; \leq \; C_8\,q^{\ell - n+1} \| \Phi_m^*\varphi_X\|_{L^2(M, \tilde{\omega} )}.
\leqno{(4.22)}
$$
In view of (4.21) and (4.22), 
$|c_X| \, 
\, \leq \, C_{9} \, q^{\ell + 2} \| \Phi_m^*\varphi_X\|_{L^2(M, \tilde{\omega} )}$ 
for some $C_{9}$.
Then by the definition of $\phi_X$ together with (4.12), we obtain
$$
\begin{cases}
&\|\phi_X\|_{L^2(M, \tilde{\omega})} \; \geq\; \| \Phi_m^*\varphi_X\|_{L^2(M, \tilde{\omega} )}
- \|c_X\|_{L^2(M, \tilde{\omega} )}\\
&\geq \; (1- C_{9} \, q^{\ell +2} \{ c_1(L)^n[M]\}^{1/2})\, \| \Phi_m^*\varphi_X\|_{L^2(M, \tilde{\omega} )}\\
& \geq \; C_{10}\| \Phi_m^*\varphi_X\|_{L^2(M, \tilde{\omega} )}
\geq  \; C_{11}\| \Phi_m^*\varphi_X\|_{L^2(M, {\omega}(\ell ) )}
\end{cases}\leqno{(4.23)}
$$
for some $C_{10}$ and $C_{11}$. Since $\theta_2 = \int_M (B_2/B_0)^2\, \omega (\ell )^n 
= m^2 \| \Phi_m^*\varphi_X\|^2_{L^2(M, {\omega}(\ell ) )}$, we see from 
(4.20) and (4.23) that
$$
\|\Phi_m^*\mathcal{X}_{TM_m^{\perp}}\|^2_{L^2(M, \tilde{\omega})}
\;\geq\; C_{12}\, q^{3}\,{\theta_2}
\leqno{(4.24)}
$$
for some $C_{12}$. Note that $q^{1/2}|\mathcal{X}_{|M_m}|^{}_{\omega_{\operatorname{FS}}}
= |\mathcal{X}_{|M_m}|^{}_{\tilde{\omega}} \geq |\mathcal{X}_{TM^{\perp}_m}|^{}_{\tilde{\omega}}$.
In view of the second line in (4.9),
it follows from (4.24)  that, for $|t| \leq \delta_0$,
$$
\begin{cases}
\;\;\ddot{f}_{X,m} (t) \;&\geq \; C_2  \int_{M_m} (q^{1/2}|\mathcal{X}_{|M_m}|_{\omega_{\operatorname{FS}}})^2
\,\omega_{\operatorname{FS}}^n\\
&\geq \; C_2\,q^{-n} \| \Phi_m^*\mathcal{X}_{TM^{\perp}_m}\|^2_{L^2(M, \tilde{\omega})}\;
\geq\; C_{13}\, q^{3-n}\,\theta_2
\end{cases} \leqno{(4.25)}
$$
for some $C_{13}$. Therefore, by (4.22) and (4.25), we obtain
$$
\dot{f}_{X,m} (\delta_0 )\; \geq \; R
\qquad \text{ and } \qquad 
\dot{f}_{X,m} (-\delta_0 )\; \geq \; -R,
$$
where $R:= C_{13}\,\delta_0 \,q^{3-n}\theta_2 - C_8 \, q^{\ell -n +1} 
\| \Phi_m^*\varphi_X\|_{L^2(M, {\omega}(\ell ) )}$.
Moreover, if $R>0$, then $\theta_2 >0$, and hence by (4.25), $\ddot{f}_{X,m}(t) >0$ for 
$|t|\leq \delta_0$.
It now suffices to show $R >0$ for $m \gg 1$. By $\theta_2  
= m^2 \| \Phi_m^*\varphi_X\|^2_{L^2(M, {\omega}(\ell ) )}$, we can write $R$ as 
$$
R \; =\; q^{3-n}\sqrt{\theta_2}\{ -\,C_8\,q^{\ell -1}  + C_{13} \,\delta_0 \sqrt{\theta_2}\}.
$$
Recall that, by our assumption of Case 2, we have $\theta_1 \leq 2 \theta_2$.
Hence by (4.17) together with the definetion of $\delta_0$, we obtain 
$$
\delta_0 \sqrt{\theta_2} \;\geq\; \delta_0 \sqrt{\theta_1}/\sqrt{2}
\;\geq\;  \sqrt{n! (1- \varepsilon)} \,q^{(n/2)+1}\; =\; C_{14}\, q^{(n/2)+1},
$$
where $C_{14} := \sqrt{n! (1- \varepsilon)}$. In view of the definition $\ell := [n/2] +3$,
$$
(n/2) + 1  \; < \; [n/2] +2\; =\;\ell -1.
$$
Therefore $R > 0$ for $m \gg 1$, as required.
\qed

\medskip\noindent
{\em Remark\/ $4.26$}: Assume that  $(M, L^{\otimes m})$ is Chow-polystable relative to 
 a maximal algebraic torus $T_{\operatorname{max}}$ in $\operatorname{Aut}^0(M)$.
 Then the arguments in \cite{MN}, Step 2, which uses \cite{Mt} allow us to obtain  finiteness 
 of the isotropy subgroup 
of $G_m$ at $[\tilde{M}_m]$. Hence in this case, $(M, L^{\otimes m})$ is 
Chow-stable relative to $T_{\operatorname{max}}$.

\section{Polybalanced metrics}

\medskip
As in the introduction, we consider a $K$-invariant Hermitian $h$ metric for $L$ such that $\omega = c_1 (L; h)$ is K\"ahler. Let $m$ be a positive integer.
In (3.2), 
we choose an orthonormal basis $\{ \sigma_{\gamma ,i}\,;\, i =  1,2,\dots, {n}_{\gamma}\}$
of the space $(\hat{V}_{m,\gamma}, \rho_m (h))$ for each $\gamma$.
In this section, we discuss the result in \cite{MM} from a slightly different point of view.
For recent related works, see \cite{ST} and \cite{H}.

\medskip\noindent
{\em Definition \/$5.1$.} For a polarized algebraic manifold $(M,L)$, $\omega$ is called an  {\it $m$-th polybalanced metric relative to $T$\/},
if for some $\mathcal{Y}\in \sqrt{-1}\,\frak t_c$ satisfying $1 - \hat{\chi}_{m,\gamma} (\mathcal{Y}) >0 $ for all 
$\gamma$, there exists
a positive real constant $C$ such that
$$
\sum_{\gamma = 1}^{\eta_m} \;\{1 - \hat{\chi}_{m,\gamma} (\mathcal{Y})\}\,
| \sigma_{\gamma ,i}|_h^2 \;=\; C.
\leqno{(5.2)}
$$ 

This concept is closely related to relative Chow-polystability.
For brevity, we use the notation in Section 4 freely until the end of this section.

\medskip\noindent
{\bf Theorem 5.3}: 
{\em $(M,L^{\otimes m})$ is Chow-polystable relative to $T$ if and only if 
$(M,L)$ admits an $m$-th polybalanced 
metric relative to $T$.}

\medskip\noindent
{\em Proof\/}: The proof of ``only if" part follows from Theorem C and (3.7) in \cite{MM}.
For ``if" part, we give a proof as follows:
Let
$$
\{\hat{\sigma}_{\gamma ,i}\,;\, i= 1,2, \dots ,n_{\gamma}, \,\gamma = 1,2,
\dots, \eta_m \}
$$ 
be the basis for $V_m$  obtained from $\{\sigma_{\gamma ,i}\,; \,i= 1,2, \dots ,n_{\gamma},\, \gamma = 1,2,
\dots, \eta_m \}$ 
by replacing $\mathcal{Y}(\ell )$ by 
$\mathcal{Y}$ in (4.2).  
Then by (4.4), we have the Kodaira embedding 
$\Phi_m : M\hookrightarrow \Bbb P^{N_m -1}(\Bbb C )$.
Also by (4.6),  we have the
function $f_{X,m}(t)$ for the orbit through $\tilde{M}_m$ of the one-parameter group
$\exp (tX)$, $t \in \Bbb R$, generated by $0\neq X \in \frak p^{\perp}$.
Since (5.2) is written as $\Sigma_{j =1}^{N_m} |\hat{\sigma}_{j} |^2_h = C$ 
in terms of the notation (4.1), 
by taking  $(\sqrt{-1}/2\pi )\partial \bar{\partial}\,{\log}$ of both sides of (5.2), we obtain
$$
\Phi_m^* \omega_{\operatorname{FS}} \, =\, m\,\omega.
$$
Put $a_{\gamma, i} := \hat{\chi}_{m,\gamma}(\mathcal{Y})$.
Then by $\mathcal{Y} \in \sqrt{-1}\, \frak t_c \subset \frak t$ and $X \in\frak p^{\perp} \subset \frak g_m
\subset \frak{sl}(V_m)$, in view of (3.1), it follows that $\Sigma_{j=1}^{N_m} b_j =\Sigma_{j=1}^{N_m} b_j a_j =0$,
where $b_j$ is as in (4.5). 
Now by replacing $h(\ell )$ by $h$ in (4.10), we obtain
\begin{align*}
\dot{f}_{X,m} (0) \;& =\; (n+1)\int_M \frac{\Sigma_{j=1}^{N_m}\, b_j |\hat{\sigma}_j |^2_{h}}{\Sigma_{j=1}^{N_m}\, |\hat{\sigma}_j |_{h}^2}\,(m\omega )^n\\
&=\; \frac{m^n(n+1)}{C}\int_M \Sigma_{j=1}^{N_m}\, b_j (1 - a_j )
|{\sigma}_j |^2_{h}\, \omega^n \\
&= \; \frac{m^n(n+1)}{C}\;\Sigma_{j=1}^{N_m}\, b_j (1 - a_j ) \; =\; 0.
\end{align*}
By this together with the convexity
 $\ddot{f}_{X, m} (t) \geq 0$, the function $f_{X,m}$ attains a minimum 
 at the origin.
Therefore every special one-parameter subgroup of $G_m$
has a closed orbit through $\tilde{M}_m$ in $W_m^*$, and hence
 the orbit $G_m\cdot \tilde{M}_m$ is closed in 
$W_m^*$ (cf.~\cite{M0}, p.568), as required.
\qed

\medskip\noindent
{\em Remark\/ $5.4$}: If $(M,L)$ is asymptotically Chow-polystable relative to $T$, then for $m \gg 1$, 
there exists an $m$-th polybalanced metric $\omega$ such that 
$$
\hat{\chi}_{m,\gamma} (\mathcal{Y} ) \; = \; O (q),
$$ 
i.e.,  the inequality $|\hat{\chi}_{m,\gamma} (\mathcal{Y} )| \leq C' q$ holds for some positive constant $C'$ 
independent of $m$, $\gamma$ and $\mathcal{Y}$ (see \cite{MM}, Theorem A).

\section{Strong relative K-stability}

In this section, for a polarized algebraic manifold $(M,L)$,
we consider an algebraic torus $T$ in $\operatorname{Aut}^0(M)$.
Let the group $ \Bbb C^*$ act on the affine line 
 $ \Bbb A^1 := \{z\in \Bbb C\}$ by multiplication of complex numbers,
$$
\Bbb C^*\times \Bbb A^1 \to \Bbb A^1,
\qquad (t, z) \mapsto
t z.
$$
By fixing a Hermitian metric $h$ for $L$ such that $\omega := c_1(L;h)$ is K\"ahler, we endow $V_{m}:= H^0(X,L^{\otimes m})$ 
with the Hermitian metric $\rho_{m}(h)$ 
as defined in the introduction.
We then consider the Kodaira embedding 
$$
 \Phi_{m} \,:\, X \, \hookrightarrow \, \Bbb P^*(V_{m}),
 \qquad x \mapsto (\tau_1 (x): \tau_2 (x): \dots : \tau_{N_{m}}(x)),
 $$
where $(\tau_1, \tau_2, \dots, \tau_{N_{m}})$ is an orthonormal basis
for $(V_{m}, \rho_{m}(h))$.
For $G_m$ as in Section 3, we consider an algebraic group homomorphism
$$
\psi \, : \,\Bbb C^* \,\to \,G_m
$$
such that the maximal compact subgroup $S^1 \subset \Bbb C^* $ acts isometrically on 
the space $(V_{m}, \rho_{m}(h))$. Put $M_m := \Phi_m (M)$.
Then by setting  
$$
\mathcal{M}^{\psi}_z := \{z\}\times\psi (z)  M_{m},
\qquad z \in \Bbb C^*,
$$
we consider the irreducible algebraic subvariety $\mathcal{M}^{\psi}$ of $\Bbb A^1 \times \Bbb P^* (V_{m})$ obtained as 
the closure of the subset
$$
\bigcup_{z\in \Bbb C^*} \;\mathcal{M}^{\psi}_z
$$ 
in $\Bbb A^1 \times \Bbb P^* (V_{m})$,
where $\psi (z)  $ in $G_m$ acts naturally on the space  
$\Bbb P^* (V_{m})$ of all hyperplanes in $V_{m}$ passing through the origin.
Let
$$
\pi : \mathcal{M}^{\psi} \to \Bbb A^1
$$ 
be the map induced by the projection of $\Bbb A^1 \times \Bbb P^* (V_{m})$ to the first factor $\Bbb A^1$.
For the hyperplane bundle $\mathcal{O}_{\Bbb P^*(V_{m})}(1)$ on $\Bbb P^*(V_{m})$, 
we consider the pullback 
$$
\mathcal{L}^{\psi}\, :=\,\operatorname{pr}_2^*\mathcal{O}_{\Bbb P^*(V_{m})}(1)_{|\mathcal{M}^{\psi}},
$$
where $\operatorname{pr}_2 : \Bbb A^1 \times \Bbb P^* (V_{m}) \to \Bbb P^* (V_{m})$
denotes the projection to the second factor.
For the dual space $V_{m}^*$ of $V_{m}$,
the $\Bbb C^*$-action on $\Bbb A^1 \times V_{m}^*$ defined by
$$
\Bbb C^* \times (\Bbb A^1 \times V_{m}^*)\to \Bbb A^1 \times V_{m}^*,
\quad (t, (z, p))\mapsto  (tz, \psi (t) p),
$$
induces $\Bbb C^*$-actions on $\Bbb A^1 \times \Bbb P^*(V_{m})$ and $\mathcal{O}_{\Bbb P^*(V_{m})}(-1)$, where $G_m$ acts on $V_{m}^*$ by the contragradient representation.
 This then induces $\Bbb C^*$-actions on $\mathcal{M}^{\psi}$ and $\mathcal{L}^{\psi}$, and hence 
 $\pi : \mathcal{M}^{\psi} \to \Bbb A^1$ is a $\Bbb C^*$-equivariant 
 projective morphism with a relatively very ample line bundle 
$\mathcal{L}^{\psi}$ satisfying
$$
(\mathcal{M}^{\psi}_z, \mathcal{L}_z^{\psi})\; \cong \; (M,L^{\otimes m}),
\qquad z \neq 0,
$$
where $\mathcal{L}_z^{\psi}$ is the restriction of $\mathcal{L}^{\psi}$ to $\mathcal{M}^{\psi}_z := \pi^{-1}(z)$.
Then a  triple $({\mathcal{M}}, {\mathcal{L}}, \psi )$ is called 
a {\it test configuration for $(M,L)$}, if we have both 
$$
\mathcal{M} = \mathcal{M}^{\psi} \quad \text{ and  }\quad  \mathcal{L}= \mathcal{L}^{\psi}.
$$
Here $m$ is called the {\it exponent} of $({\mathcal{M}}, {\mathcal{L}}, \psi )$.
A test configuration $({\mathcal{M}}, {\mathcal{L}}, \psi )$ is called {\it trivial},
if $\psi$
is a trivial homomorphism.
Let {\bf M} be the set of all sequences $\{\mu_j \}$ of test configurations 
$$
\mu_j  \, =\, (\mathcal{M}_j, \mathcal{L}_j, \psi_j ), \qquad j =1,2,\dots,
$$
for $(M,L)$  
such that the exponent $m_j$ of the test configuration $\mu_j$ satisfies the following 
growth condition:
$$
\text{$m_j \to +\infty$,  \;\; as $j \to \infty$.}
$$
In \cite{M4}, to each $\{\mu_j\}\in  $
{\bf M}, we associated the Donaldson-Futaki invariant 
$$
F_1 (\{\mu_j\}) \in  \Bbb R \cup \{-\infty\},
$$
which is viewed as a generalization
of the Donaldson-Futaki invariant $DF(\mu )$ of a test configuration $\mu$.
We can also define the following strong version of 
K-stability and K-semistability:

\medskip\noindent
{\em Definition\/ $5.1$}. 
(1) A polarized algebraic manifold $(M,L)$ is called 
{\it strongly $K$-semistable relative to $T$}, 
if $F_1(\{\mu_j\} ) \leq 0$ for all 
$\{\mu_j\} \in $
{\bf M}.

\medskip\noindent
(2) A strongly K-semistable polarized algebraic manifold $(M,L)$
 is called 
{\it strongly $K$-stable relative to $T$},
if  for every $\{\mu_j\}\in $ {\bf M} with $ F_1(\{\mu_j\} ) = 0$,
there exists $j_0$ such that
 $\mu_j$ are trivial for all $j$ with $j \geq j_0$.
 
 \section{Concluding Remarks}
 
 The Yau-Tian-Donaldson Conjecture for K\"ahler-Einstein cases was solved affirmatively by 
 Chen-Donaldson-Sun \cite{CD} and Tian \cite{T1}. However, for general polarization cases or extremal 
 K\"ahler cases, the conjecture is still open. In this paper, we discuss extremal  K\"ahler versions of this conjecture by focussing on
 the difference between strong K-stability and K-stability.
 For an arbitrary polarized algebraic manifold $(M,L)$ as in the introduction, recall the following definition of 
 K-stability 
  \cite{D1}  (cf.~\cite{T0}) :
 
 \medskip\noindent
{\em Definition\/ $6.1$}. (1) $(M,L)$ is called 
{\it $K$-semistable}, 
if $DF(\mu ) \leq 0$  for all 
test configurations $\mu$ for $(M,L)$.

\smallskip\noindent
(2) A K-semistable $(M,L)$
 is called 
{\it $K$-stable},
if  every test configuration $\mu$ for $(M,L)$ with $ DF(\mu ) = 0$
is trivial.

\medskip
We now consider a maximal algebraic torus $T_{\operatorname{max}}$ in $\operatorname{Aut}^0(M)$.
 As to the existence of extremal K\"ahler metrics, 
it is natural to ask the following: 

\medskip\noindent
{\bf Conjecture I}: {\em A polarized algebraic manifold $(M,L)$ is K-stable relative to $T_{\operatorname{max}}$
if and only if $c_1(L)$ admits an extremal K\"ahler metric.}

\medskip\noindent
{\bf Conjecture II}: {\em A polarized algebraic manifold $(M,L)$ is strongly K-stable relative to $T_{\operatorname{max}}$
if and only if $c_1(L)$ admits an extremal K\"ahler metric.}

\medskip
By a result of Donaldson \cite{D}, every constant scalar curvature K\"ahler metric is approximated by a sequence of 
balanced metrics. In other words, a balanced metric can be viewed as a quantized version of a constant scalar curvature 
K\"ahler metric. Similarly,  a polybalanced metric 
can be viewed as  a quantized version of an extremal K\"ahler metric. Since  
existence of polybalanced metrics 
corresponds to 
relative Chow-polystability, the following  fact is viewed as
a quantized version of the existence part of Conjecture II:

\medskip\noindent
{\bf Fact}: (cf.~\cite{MN}) {\em  If a polarized algebraic manifold $(M,L)$ is strongly K-stable relative to an algebraic 
torus $T$ in $\operatorname{Aut}^0(M)$, then
$(M, L^{\otimes m})$, $m \gg 1$,  are Chow-stable relative to $T$.}

\medskip
Moreover, we expect that ``if" part  of Conjecture II is true. This will be discussed in a forthcoming paper \cite{M5}.
In view of the above Fact, by assuming  that ``if" part of Conjecture II is true, we immediately obtain the case $T = T_{\operatorname{max}}$ 
of Corollary in the introduction.

\medskip
For a polarized algebraic manifold $(M,L)$, 
let  $T_{\operatorname{ex}}$ be the algebraic torus   in $\operatorname{Aut}^0(M)$ generated by the extremal K\"ahler vector field.
By a recent result of Yotsutani-Zhou \cite{YZ}, a smooth polarized toric Fano threefold
$$
\Pi := (\mathcal{E}_4, K^{-1}_{\mathcal{E}_4})
$$
is K-stable relative to $T_{\operatorname{max}}$, and is not asymptotically Chow-stable 
relative to $T_{\operatorname{ex}}$. 
Let $Z_{\Bbb  C}$ be as in the introduction.
We finally pose the following:

\medskip\noindent
{\bf Problem}: {\em Check whether or not $\Pi$ is asymptotically Chow stable relative to $Z_{\Bbb  C}$.
Or more generally, clarify whether there is an example of a polarized algebraic manifold $(M,L)$ which is K-stable relative to $T_{\operatorname{max}}$ and is not asymptotically Chow-stable relative to $Z_{\Bbb C}$.
}

\medskip
 If there is such an example of a polarized algebraic manifold $(M,L)$, then by Corollary in the introduction, $c_1(L)$ admits no extremal K\"ahler metrics. In other words, this gives a counter-example to Conjecture I above.

\bigskip\noindent
{\footnotesize
{\sc Mathematics Department}, {\sc Osaka University}, {\sc Toyonaka, Osaka, 
560-0043 Japan}}


\begin{thebibliography}{Mat}

\bibitem{A}
{\sc V.~Apostolov and H.~Huang}:
{\it A splitting theorem for extremal K\"ahler metrics},
J. Geom. Anal. {\bf 25} (2015), 149--170.

\bibitem{C}
{\sc E.~Calabi}: 
{\it Extremal K\"ahler metrics II}, in ``Differential Geometry and Complex Analysis" (ed. I.~Chavel, 
H.M.~Farkas), Springer-Verlag, 1985, 95--114.

\bibitem{Ca}
{\sc D.~Catlin}:
{\it The Bergman kernel and a theorem of Tian},
in ``Analysis and Geometry in Several Complex Variables'' (ed. G.~Komatsu, M.~Kuranishi), 
Trends in Math., Birkh\"auser, 1999, 1--23.

\bibitem{CD}
{\sc X.~Chen, S.K.~Donaldson and S.~Sun}:
{\it K\"ahler-Einstein metrics on Fano manifolds, I: Approximation of metrics with cone singularities, 
II: Limits with cone angle less than $2\pi$, III: Limits as cone angle approaches $2\pi$ and completion of the main proof}, J.~Amer.~Math.~Soc. {\bf 28} (2015), 183--197, 199--234, 235--278.

\bibitem{D}
{\sc S.K.~Donaldson}:
{\it Scalar curvature and projective embeddings, I}, \,J. Differential Geom.
{\bf 59} (2001), 479--522.

\bibitem{D1}
{\sc S.K.~Donaldson}:
{\it Scalar curvature and stability of toric varieties},  J. Differential Geom.
{\bf 62} (2002), 289--349.

\bibitem{F}
{\sc A.~Futaki}:
{\it Asymptotic Chow semi-stability and integral invariants},
Internat. J. Math. {\bf 15} (2004), 967--979.


\bibitem{FM} 
{\sc A.~Futaki and T.~Mabuchi}:
{\it Bilinear forms and extremal K\"ahler vector fields associated with K\"ahler classes},
Math. Ann. {\bf 301} (1995), 199--210.

\bibitem{H}
{\sc Y.~Hashimoto}:
{\it Quantisation of extremal K\"ahler metrics},
arXiv: 1508.02643, math.DG (2015).

\bibitem{L}
{\sc Z.~Lu}: 
{\it On the lower order terms of the asymptotic expansion Tian-Yau-Zelditch},
Amer. J. Math. {\bf 122} (2000), 235--273.

\bibitem{MA}
{\sc T.~Mabuchi}:
{\it An algebraic character associated with Poisson brackets},
in ``Recent Topics in Differential and Analytic Geometry'' (ed. T.~Ochiai),
Adv. Stud. Pure Math. {\bf 18-I}, Kinokuniya and Academic Press, 1990, 339--358.

\bibitem{M0}
{\sc T.~Mabuchi}:
{\it Stability of extremal K\"ahler manifolds},
Osaka J.~Math. {\bf 41} (2004), 563--582.

\bibitem{M1}
{\sc T.~Mabuchi}:
{\it An energy-theoretic approach to the Hitchin-Kobayashi correspondence for manifolds, II},
\,arXiv: 0410239, math.DG (2004).

\bibitem{M2}
{\sc T.~Mabuchi}:
{\it An energy-theoretic approach to the Hitchin-Kobayashi correspondence for manifolds, I},
 \,Invent.~Math. {\bf 159} (2005), 225--243.

\bibitem{M3}
{\sc T.~Mabuchi}:
{\it An energy-theoretic approach to the Hitchin-Kobayashi correspondence for manifolds, II}, 
\,Osaka J.~Math. {\bf 46} (2009), 115--139.

\bibitem{MM}
{\sc T.~Mabuchi}:
{\it Asymptotics of polybalanced metrics under relative stability constraints}, 
\,Osaka J.~Math. {\bf 48} (2011), 845--856.


\bibitem{M4}
{\sc T.~Mabuchi}:
{\it The Donaldson-Futaki invariant for sequences of test configurations},
in ``Geometry and Analysis on Manifolds," 
 \,Progr.~Math. {\bf 308}, Birkh\"auser Boston, 2015, 395--403.
 
 \bibitem{M5}
 {\sc T.~Mabuchi}:
 {\it A stronger concept of K-stability},
 A revised version of arXiv: 0910.4617, math.DG (2009),
 in preparation.
 
 
 \bibitem{MN}
 {\sc T.~Mabuchi and Y.~Nitta}:
 {\it Strong K-stability and asymptotic Chow-stability},
 in ``Geometry and Analysis on Manifolds,"
 \,Progr.~Math. {\bf 308}, Birkh\"auser Boston, 2015, 405--411.
 
 \bibitem{Mt}
 {\sc Y.~Matsushima}:
{\it  Espaces homog\`enes de Stein des des groupes de Lie complexes}, 
Nagoya~Math.~J. {\bf 18} (1961), 153--164.
 
\bibitem{OSY}
{\sc H.~Ono, Y.~Sano and N.~Yotsutani}:
{\it  An example of an asymptotically Chow unstable manifold with constant scalar curvature},
Ann. Inst. Fourier (Grenoble) {\bf 62} (2012), 1265--1287.

\bibitem{PS}
{\sc D.H.~Phong and J.~Sturm}:
{\it Scalar curvature, moment maps, and the Deligne pairing},
\,arXiv: 0209098, math. DG (2002).

\bibitem{ST}
{\sc Y.~Sano and C.~Tipler}:
{\it Extremal metrics and lower bound the modified K-energy},
J. Eur. Math. Soc. {\bf 17} (2015), 2289--2310.

\bibitem{S1}
{\sc G.~Sz\'ekelyhidi}:
{\it Extremal metrics and K-stability}, Dissertation, Imperial College, London, arXiv: 0611002, math.DG (2006).

\bibitem{S2}
{\sc G.~Sz\'ekelyhidi}:
{\it Extremal metrics and K-stability}, Bull.~London Math.~Soc. {\bf 39} (2007), 76--84.

\bibitem{T}
{\sc G.~Tian}:
{\it On a set of polarized K\"ahler metrics on algebraic manifolds}, J.~Differential~Geom. {\bf 32} (1990), 99--130.

\bibitem{T0}
{\sc G.~Tian}:
{\it K\"ahler-Einstein metrics with positive scalar curvature}, Invent.~Math. {\bf 130} (1997), 1--37.

\bibitem{T1}
{\sc G.~Tian}:
{\it K-stability and K\"hler-Einstein metrics}, Comm.~Pure~Appl.~Math. {\bf 68} (2015), 1085--1156,
 2082--2083.
 
 \bibitem{YZ}
 {\sc N.~Yotsutani and B.~Zhou}:
 {\it Relative algebro-geometric stabilities of toric manifolds}, \, arXiv: 1602.08201, math. DG (2016).

\bibitem{Ze}
{\sc S.~Zelditch}:
{\it Szeg\"o kernels and a theorem of Tian},
Internat. Math. Res. Notices {\bf 6} (1998), 317--331.

\bibitem{Z}
{\sc S.~Zhang}: Heights and reductions of semi-stable varieties, Compos.~Math. {\bf 104} (1996), 77--105.

\end{thebibliography}
\end{document}